\documentclass[a4paper]{article}

\usepackage[plainpages=false, colorlinks=true, 
            linkcolor=black, urlcolor=black, citecolor=black]{hyperref}
\usepackage{geometry}
\usepackage{tabularx}
\usepackage{tikz}
\usepackage{overpic}
\usepackage{amsmath, amssymb, amsthm,  amsfonts, mathrsfs, cite}
\usepackage{multirow}
\usepackage{color}
\usepackage{booktabs}
\usepackage{rotating}
\usepackage{tikz}
\usepackage{soul}
\usepackage{bm}
\usepackage{hyperref}
\usepackage{scrextend}
\usepackage{arydshln}
\usepackage {url}
\usepackage[font=scriptsize, labelfont=bf]{caption}

\usepackage{graphicx} 
\usepackage{epstopdf} 
\usepackage{algpseudocode}
\usepackage{graphicx}
\usepackage{array,  tabularx}
\usepackage{algorithm, algcompatible}
\usepackage{booktabs} 
\usepackage{paralist} 
\usepackage{verbatim} 
\usepackage{subfig} 
\usepackage{xcolor, marginnote, enumitem}
\usepackage[numbers]{natbib}
\numberwithin{equation}{section}
\theoremstyle{theorem}
\newtheorem{lemma}{Lemma}
\newtheorem{theorem}{Theorem}
\newtheorem{proposition}{Proposition}

\theoremstyle{remark}
\newtheorem{remark}{Remark}

\theoremstyle{definition}
\newtheorem{definition}{Definition}

\DeclareMathOperator{\dist}{dist}

\newcommand{\tbfzero}{\textbf{0}}

\newcommand{\lookUp}[1]{}

\newcommand{\bitem}{\begin{itemize}}
	\newcommand{\eitem}{\end{itemize}}

\newcommand{\bpm}{\begin{pmatrix}}

\DeclareMathAlphabet{\mathbfit}{OML}{cmm}{b}{it}

\usepackage{xspace}
\usepackage{bold-extra}
\usepackage[most]{tcolorbox}

\colorlet{texcscolor}{blue!50!black}
\colorlet{texemcolor}{red!70!black}
\colorlet{texpreamble}{red!70!black}
\colorlet{codebackground}{black!25!white!25}

\date{}

\begin{document}

\title{A preconditioned third-order implicit-explicit algorithm with a difference of varying convex functions and extrapolation}

\author{Kelin Wu\thanks{School of Mathematics, 
		Renmin University of China,  China  \  \href{mailto:kelinwu@ruc.edu.cn}{kelinwu@ruc.edu.cn}}
    \and   Hongpeng Sun\thanks{School of Mathematics,
Renmin University of China, China \  \href{mailto:hpsun@amss.ac.cn}{hpsun@amss.ac.cn} }
}

\maketitle

\begin{abstract}
This paper proposes a novel preconditioned implicit-explicit algorithm enhanced with the extrapolation technique for nonconvex optimization problems. The algorithm employs a third-order Adams-Bashforth scheme for the nonlinear and explicit parts and a third-order backward differentiation formula for the implicit part of the gradient flow in variational functions. The proposed algorithm, akin to a generalized difference-of-convex approach, employs a changing set of convex functions in each iteration. Under the Kurdyka-\L ojasiewicz  properties, the global convergence of the algorithm is guaranteed, ensuring that it converges within a finite number of preconditioned iterations. Our numerical experiments, including least squares problems with SCAD regularization and the graphical Ginzburg-Landau model, demonstrate the proposed algorithm's highly efficient performance compared to conventional difference-of-convex algorithms.
\end{abstract}

\paragraph{MSCodes.}{
  65K10,  
65F08,  
49K35, 
90C26 
}


\section{Introduction}\label{sec:intro}
We focus on the following nonconvex optimization problem
\begin{equation}\label{eq:basefunctional}
    \min_{x \in X} E(x)= H(x) + F(x)
\end{equation}
where $H(x)$ is a proper closed convex function and $F(x)$ is an $L$-smooth function, i.e. $\nabla F(x)$ is Lipschitz continuous with constant $L$. $X$ is a finite-dimensional Hilbert space. nonconvex optimization problems are ubiquitous in modern computational mathematics, with critical applications spanning machine learning, signal processing, finance, and beyond \cite{APX,Pangcui,BF}. The difference of convex functions algorithm (DCA) is a powerful method for tackling this class of problems by leveraging the problem's structure \cite{APX,Shensun2023,2020boosted}. 

This work is inspired by the implicit-explicit (IMEX) approach of treating stiff dissipative terms implicitly and non-stiff terms explicitly~\cite{ARW1995}. The widely adopted third-order IMEX schemes are of vital significance for the phase field crystal model and demonstrate robust applicability in simulating crystal growth, microstructure evolution~\cite{newthirdorder,thirdorderenergy,thirdorderSD}, and convection-diffusion problems~\cite{ARW1995}. The method can be written as \cite{ARW1995,Feng2013}
\begin{equation}\label{eq:3ord:ori}
\frac{1}{k} \left( \frac{11}{6} x^{n+1} - 3x^n + \frac{3}{2} x^{n-1} - \frac{1}{3} x^{n-2} \right) = -3f(x^n) + 3f(x^{n-1}) - f(x^{n-2}) -\nabla H(x^{n+1}), 
\end{equation}
where $f(x)=\nabla F(x)$, $\nabla H$ is continuous, $k$ is a positive step size. The explicit scheme is 
a third-order Adams-Bashforth scheme for the nonlinear and explicit parts of the gradient flow, while the implicit scheme is a third-order backward differentiation formula for the implicit part of the gradient flow.

With a scaling to the time step, we propose the following scheme with a nonsmooth function $H$:
\begin{equation}\label{eq:3ord:nonsmooth}
    0\in\frac{12}{11\delta t}\left[\frac{11}{6}x^{n+1}-3x^n+\frac{3}{2}x^{n-1}-\frac{1}{3}x^{n-2}\right]+\partial H(x^{n+1})+[3f(x^n)-3f(x^{n-1})+f(x^{n-2})],
\end{equation}
where $\partial H(x)$ is the subdifferential of $H(x)$ \cite{HBPL,KK} and $\delta t>0$ refers to a step size. Motivated by the refined difference-of-convex (DC) approach with extrapolation \cite{Wen2018}, as well as the preconditioned framework for DCA and second-order splitting algorithm \cite{DS,shen2024secondorder}, we proposed a preconditioned third-order implicit-explicit algorithm with extrapolation to solve \eqref{eq:basefunctional}. For subsequent analysis, we introduce the following energies:
\begin{equation}\label{eq:def:en:hn}
\begin{aligned}
     E^n(x)&=H^n(x)-F^n(x),\\
    H^n(x)&=H(x)+\frac{1}{\delta t}\|x-x^n\|^2,\\
    F^n(x)&=\frac{\alpha}{\delta t}\|x-x^n\|^2+\frac{14}{11\delta t}\langle x-x^n,x^n-x^{n-1}\rangle-\frac{4}{11\delta t}\langle x-x^{n-1},x^{n-1}-x^{n-2}\rangle-F(x)\notag\\
    &\quad -2\langle f(x^n)-f(x^{n-1}),x-x^{n-1}\rangle+\langle f(x^{n-1})-f(x^{n-2}),x-x^{n-2}\rangle,
\end{aligned}
\end{equation}
where $\alpha$ is a positive constant.

Now we incorporate extrapolation techniques to accelerate the algorithm. Set 
\begin{align}
    y^n=x^n+\beta_n(x^n-x^{n-1})
\end{align}
where
\begin{equation}\label{beta}
    \{\beta_n\}_n \subset [0,1),\quad \sup_n \beta_n < 1
\end{equation}
before using a proximal term to update $x^{n+1}$:
\begin{align} \label{eq:mini:dc:sub}
    x^{n+1}=\underset{x}{\arg\min} \left\{H^n(x)-\langle \nabla F^n(x^n),x\rangle+\frac{1}{2}\|x-y^n\|_M^2\right\}.
\end{align}
Here, the positive semi-definite and symmetric weight $M$ is utilized to create efficient preconditioners and $\|x\|_{M}^2:=\langle Mx, x\rangle$.  By applying the first-order optimality condition, the minimization problem for determining $x^{n+1}$ leads to the equation:
\begin{equation}\label{eq:solve:y:pre:general}
\begin{aligned}
        0\in&\frac{12}{11\delta t}\left[\frac{11}{6}x^{n+1}-3x^n+\frac{3}{2}x^{n-1}-\frac{1}{3}x^{n-2}\right]+M(x^{n+1}-y^n)\\
    &+\partial H(x^{n+1})+[3f(x^n)-3f(x^{n-1})+f(x^{n-2})].
\end{aligned}
\end{equation}
Equivalently, for $w^{n+1}\in\partial H(x^{n+1})$, we have
\begin{equation}\label{eq:solve:y:pre}
\begin{aligned}
        &\frac{12}{11\delta t}\left[\frac{11}{6}x^{n+1}-3x^n+\frac{3}{2}x^{n-1}-\frac{1}{3}x^{n-2}\right]+M(x^{n+1}-y^n)+w^{n+1}\\
        &+[3f(x^n)-3f(x^{n-1})+f(x^{n-2})]=0.
\end{aligned}
\end{equation}

It is worth pointing out that the conditions \eqref{beta} on the sequence $\{\beta_n\}_n$ are sufficiently general to accommodate a variety of extrapolation parameters. Specifically, in these schemes, the initial values are set as $\theta_{-1} = \theta_0 = 1$. For $n \geq 0$, the sequence $\{\beta_n\}_n$ is recursively defined by
\begin{equation} \label{extrapara}
    \beta_n = \frac{\theta_{n-1} - 1}{\theta_n} \quad \text{with} \quad \theta_{n+1} = \frac{1 + \sqrt{1 + 4\theta_n^2}}{2}.
\end{equation}
We also claim that the values of $\theta_{n-1}$ and $\theta_n$ should be reset to 1 under suitable conditions. In this paper, the reset $\theta_{n-1} = \theta_n = 1$ is triggered whenever the condition 
\begin{equation} \label{restart}
  \langle y^{n} - x^n, x^n - x^{n-1} \rangle > 0
\end{equation}
is satisfied. The choice of the sequence $\{\beta_n\}_n$ will be employed in our numerical experiments.

Specifically, when the functions $f$ and $g$ in the algorithm correspond to the gradients of certain potential functions, this framework can be regarded as a discrete implementation of high-order nonlinear gradient flow systems \cite{ARW1995}. 

We summarize our contributions as follows. Above all, we are the first to establish the global convergence of the iteration sequence of the third-order implicit-explicit method with function $h$ being nonsmooth. While third-order implicit-explicit methods have found widespread application in phase field crystal modeling, energy-stable simulations, and mass-conserving systems \cite{newthirdorder,thirdorderSD,thirdorderenergy}, existing convergence analyses have mainly concentrated on stability of the energy or perturbed energy \cite[Theorem 3.2]{Feng2013}, \cite[Lemma 2.3]{Shen2010}, or \cite[Theorem 1.2]{Li2017}. Advancements in Kurdyka-\L ojasiewicz (KL) analysis  \cite{Attouch2009, ABS, Li2018} have enabled us to demonstrate global convergence of the iteration sequence under mild conditions. Additionally, we have introduced extrapolation acceleration \cite{Wen2018,refinedextra,yang2024proximal,boct2022extrapolated} to enhance the third-order implicit-explicit method and incorporated a preconditioning technique \cite{DS, Shensun2023} to handle large-scale linear subproblems efficiently during each iteration. By conducting a finite number of preconditioned iterations without error control, we can guarantee global convergence through extrapolation accelerations. Numerical experiments,  including data classification and image segmentation,  have demonstrated the significant efficiency of our proposed preconditioned third-order implicit-explicit with extrapolation.

The remainder of this paper is organized as follows. Section \ref{sec2} is designed for a theoretical framework, presenting our preconditioned third-order implicit-explicit method with extrapolation and establishing its global convergence through KL properties. Section \ref{sec:pre} demonstrates the algorithm's effectiveness through numerical experiments. Finally, Section \ref{sec:conclusion} gives a conclusion.

\section{Convergence analysis of the proposed algorithm}\label{sec2}
\subsection{Preliminary}
For the convergence analysis, we need the following  Kurdyka-\L ojasiewicz (KL) property and KL exponent. The KL properties facilitate the global convergence of iterative sequences, while the KL exponent aids in determining a local convergence rate. 
\begin{definition}[KL property, KL function {\cite[Definition 2.4]{ABS}} {\cite[Definition 1]{Artacho2018}} and KL exponent {\cite[Remark 6]{Bolte2014}}]\label{def:KL}Let $h: \mathbb{R}^n \rightarrow \mathbb{R}$ be a closed proper function. $h$ is said to satisfy the KL property if for any critical point $\bar x$,  there exists $\nu \in (0,+\infty]$, a neighborhood $\mathcal{O}$ of $\bar x$, and a continuous concave function $\psi: [0,\nu) \rightarrow [0,+\infty)$ with $\psi(0)=0$ such that:
	\begin{itemize}
		\item [{\rm{(i)}}] $\psi$ is continuous differentiable on $(0,\nu)$ with $\psi'>0$ over $(0,\nu)$;
		\item [{\rm{(ii)}}] for any $x \in \mathcal{O}$ with $h(\bar x) < h(x) <h(\bar x) + \nu$, one has
		\begin{equation}\label{eq:kl:def}
		\psi'(h(x)-h(\bar x)) \cdot\text{dist}(\textbf{0},\partial h(x))\geq 1 . 
		\end{equation}
	\end{itemize}
	Furthermore, for a differentiable function $h$ satisfying the KL property,  if $\psi$ in \eqref{eq:kl:def} can be chosen as $\psi(s) = cs^{1-\theta}$ for some $\theta \in [0,1)$ and $c>0$, i.e., there exist  $\bar c, \epsilon >0$ such that
	\begin{equation}\label{eq:KL:exponent:theta:exam}
	\text{dist}(\textbf{0},\partial h(x)) \geq \bar c (h(x)-h(\bar x))^{\theta}
	\end{equation}
 whenever $\|x -\bar x\| \leq \epsilon$ and $h(\bar x) < h(x) <h(\bar x) + \nu$, then we say that $h$ has the KL property at $\bar x$ with exponent $\theta$.
\end{definition}
The KL exponent is determined exclusively by the critical points. If $h$ exhibits the KL property with an exponent $\theta$ at any critical point $\bar x$, then $h$ is a KL function with an exponent $\theta$ at all points in dom $\partial h$ \cite[Lemma 2.1]{Li2018}.

It is also assumed that the energy $E(x)$ is level-bounded. This means that for every scalar $\alpha \in \mathbb{R}$, $\text{lev}_{\leq \alpha}(E): = \{x: E(x) \leq \alpha\}$ is bounded (or possibly empty) \cite[Definition 1.8]{Roc1}. If a function $f: \mathbb{R}^n \rightarrow \mathbb{R}$ is coercive (i.e., $f(x) \rightarrow +\infty$ as $\|x\|\rightarrow +\infty$), it is also level-bounded.

We now direct our attention to the preconditioning technique for the case where \( w^n =\nabla H(x^n)= A x^n - b_0 \), where \( A \) is a bounded, linear, positive semi-definite operator and \( b_0 \in X \) is a known quantity. The central idea involves employing established preconditioning techniques, such as the symmetric Gauss–Seidel, Jacobi, and Richardson methods, to solve large-scale linear systems efficiently. A key observation, which will be rigorously demonstrated, is that any finite and feasible preconditioned iteration is sufficient to guarantee global convergence of the overall third-order implicit-explicit method.

More discussion on preconditioning techniques for nonlinear convex problems can be found in \cite{BSCC}. For additional perspectives on DCA, please refer to \cite{DS, Shensun2023}. To explain the fundamental idea of preconditioning techniques, we present the classical preconditioning technique via the following proposition.

\begin{proposition}\label{prop:proximal_to_pre}

 If $H$ is quadratic function with $\nabla H|_{x=x^n}=Ax^n-b_0$, the equation \eqref{eq:solve:y:pre} can be transformed as the following preconditioned iteration
    \begin{align}
       x^{n+1}=y^n+\mathbb{M}^{-1}(b^n-Ty^n)
    \end{align}
where 
\begin{align*}
    &T=\frac{2}{\delta t}I+A, \quad  \mathbb{M}=T+M,\\
    &b^n=b_0+\frac{12}{11\delta t}(3x^n-\frac{3}{2}x^{n-1}+\frac{1}{3}x^{n-2})-[3f(x^n)-3f(x^{n-1})+f(x^{n-2})].  
\end{align*}

\end{proposition}
\begin{proof}
    With \eqref{eq:solve:y:pre:general}, we have 
    \begin{align}
        &(\frac{2}{\delta t}I+A)x^{n+1}+Mx^{n+1} \notag \\
        =&My^n+b_0+\frac{12}{11\delta t}\left(3x^n-\frac{3}{2}x^{n-1}+\frac{1}{3}x^{n-2}\right)-[3f(x^n)-3f(x^{n-1})+f(x^{n-2})].
    \end{align}
    With our notations, the equality is equivalent to
    \begin{align*}
        (T+M)x^{n+1}=My^n+b^n
        \Leftrightarrow & \mathbb{M}x^{n+1}=\mathbb{M}y^n+b^n-Ty^n\\
        \Leftrightarrow &x^{n+1}=y^n+\mathbb{M}^{-1}(b^n-Ty^n)
    \end{align*}
    which leads to the proposition.
\end{proof}

It is important to note that classical preconditioners, such as symmetric Gauss-Seidel (SGS) preconditioners, do not require the explicit specification of $M$ as in \eqref{eq:mini:dc:sub} \cite{BSCC}. The positive semi-definiteness of $M$ is inherently fulfilled through SGS iteration in solving the linear equation $Tx = b^n$. 

Let us turn to the constraint of the step size $\delta t$. We mainly need the convexity of $F^n$ leading to the following lemma. 
\begin{lemma} \label{lem:con:L}
     For the strong convexity of $F^n(x)$, supposing that the Lipschitz constant of $f(x)$ $L>0$, we need $\delta t<\frac{2\alpha}{L}$, where $\alpha$ is the positive constant in the definition of $F^n(x)$.
\end{lemma}
    \begin{proof}
     With the form of $F^n(x)$, we derive
    \begin{align*}
        \nabla F^n(x)=&\frac{2\alpha}{\delta t}(x-x^n)+\frac{14}{11\delta t}(x^n-x^{n-1})-\frac{4}{11\delta t}(x^{n-1}-x^{n-2})-f(x)\\
        &-[2f(x^n)-3f(x^{n-1})-f(x^{n-2})].
    \end{align*}
    It can be checked that for any $x_1,x_2\in X$, we have
    \begin{align*}
       &\langle \nabla F^n(x_1)-\nabla F^n(x_2),x_1-x_2\rangle
       =\langle\frac{2\alpha}{\delta t}(x_1-x_2)-f(x_1)+f(x_2),x_1-x_2\rangle\\
       &=\frac{2\alpha}{\delta t}\|x_1-x_2\|^2-\langle f(x_1)-f(x_2),x_1-x_2\rangle
       \geq (\frac{2\alpha}{\delta t}-L)\|x_1-x_2\|^2.
    \end{align*}
    Consequently, $F^n$ is strongly monotone with parameter $\frac{2\alpha}{\delta t}-L$ under the condition $\delta t<\frac{2\alpha}{L}$.
    \end{proof}
    We now proceed to present the algorithm. It is worth noting that the step size $\delta t$ is of paramount importance. Actually, along with the necessity of $\delta t<\frac{2\alpha}{L}$ to ensure the strong convexity of $F^n$,  the global convergence of our algorithm to be discussed later will add more constraints on $\delta t$, i.e, $0<\delta t<\frac{8}{77L}$. We refer to Lemma \ref{lem:square:summable} in the next subsection for details.  With these preparations, we can present the following Algorithm \ref{alg:cap} with the constraint $0<\delta t<\frac{8}{77L}$ and $\alpha>\frac{4}{77}$ for completeness.

\begin{algorithm} 
\renewcommand{\algorithmicrequire}{\textbf{Input:}}
\renewcommand{\algorithmicensure}{\textbf{Output:}}
\caption{Algorithmic framework for third-order BDF and Adams-Bashforth implicit-explicit method with extrapolation and preconditioning (shortened as \protect{\textbf{3BapDCA$_{\text{e}}$}})}\label{alg:cap} 
\begin{algorithmic}[1] 
\State Choose $x^0$, $\alpha>4/77$,  $0<\delta t<8/(77L)$. Set $x^{-1}=x^1=x^0$, $\theta_{-1} = \theta_0 = 1$.
\State $H^n(x) := H(x) + \frac{1}{\delta t}\|x-x^n\|^2$
\State $F^n(x) := \frac{\alpha}{\delta t}\|x-x^n\|^2+\frac{14}{11\delta t}\langle x-x^n,x^n-x^{n-1}\rangle-\frac{4}{11\delta t}\langle x-x^{n-1},x^{n-1}-x^{n-2}\rangle-F(x) -2\langle f(x^n)-f(x^{n-1}),x-x^{n-1}\rangle+\langle f(x^{n-1})-f(x^{n-2}),x-x^{n-2}\rangle$
\State  Define the sequence $\{\beta_n\}_n$ recursively by \eqref{extrapara} and set $y^n=x^n+\beta_n(x^n-x^{n-1})$.
\State Solve the following subproblem 
\begin{equation}\label{eq:algorithmic:update}
     x^{n+1}=\underset{x}{\arg\min} \left\{H^n(x)-\langle \nabla F^n(x^n),x\rangle+\frac{1}{2}\|x-y^n\|_M^2\right\}.
\end{equation}
The update of $x^{n+1}$ becomes the preconditioned iteration $x^{n+1} = y^n + \mathbb{M}^{-1}(b^n - Ty^n)$ as in Proposition \ref{prop:proximal_to_pre} while $w^n = \nabla H|_{x=x^n}= Ax^n-b_0$.
\State If $\langle y^{n} - x^n, x^n - x^{n-1} \rangle > 0$, reset $\theta_{n-1} = \theta_n = 1$ and go to Step 4.
\State If any given stopping criterion is satisfied, STOP and RETURN $x^{n}$; otherwise, set $n = n+1$ and go to Step 4.
\end{algorithmic}
\end{algorithm}

We then turn to the discussion of the global convergence of Algorithm \ref{alg:cap} in the next subsection.

\subsection{Global convergence}
We begin our convergence analysis with the following lemma, illustrating the relationship of the energy function between two consecutive steps.
\begin{proposition}\label{prop:descent}
Let $\{x^n\}_n$ be the sequence generated by solving \eqref{eq:algorithmic:update} and $w^n \in \partial H(x^n)$. The following inequality holds
\begin{equation}\label{lowerbound}
\begin{aligned}
    E^n(x^{n+1})\leq &E^n(x^n) - \left(\frac{\alpha+1}{\delta t}-\frac{L}{2}\right)\|x^{n+1}-x^n\|^2-(1-\frac{\beta_n}{2})\|x^{n+1}-x^n\|_M^2\\
    &+\frac{\beta_n}{2}\|x^n-x^{n-1}\|_M^2.  
\end{aligned}
\end{equation}
\end{proposition}
\begin{proof}
   Due to the strong convexity of $H^n(x)$ and $F^n(x)$, we get
    \begin{align*}
        H^n(x^n)-H^n(x^{n+1})&\geq \langle w^{n+1}+\frac{2}{\delta t}(x^{n+1}-x^n),x^n-x^{n+1}\rangle+\frac{1}{\delta t}\|x^n-x^{n+1}\|^2\\
        F^n(x^{n+1})-F^n(x^n)&\geq \langle \nabla F^n(x^n),x^{n+1}-x^n\rangle+\left(\frac{\alpha}{\delta t}-\frac{L}{2}\right)\|x^{n+1}-x^n\|^2.
    \end{align*}
     Summing the inequalities above and using the first-order optimality condition
    \begin{align*}
        w^{n+1}+\frac{2}{\delta t}(x^{n+1}-x^n)-\nabla F^n(x^n)+M(x^{n+1}-y^n) =0,
    \end{align*}
we have
    \begin{align*}
        &E^n(x^n)-E^n(x^{n+1}) \geq \left(\frac{\alpha+1}{\delta t}-\frac{L}{2}\right)\|x^{n+1}-x^n\|^2+\langle M(x^{n+1}-x^n),x^{n+1}-y^n\rangle\\
        &=\left(\frac{\alpha+1}{\delta t}-\frac{L}{2}\right)\|x^{n+1}-x^n\|^2+\langle M(x^{n+1}-x^n),x^{n+1}-x^n-\beta_n(x^n-x^{n-1})\rangle\\
        &=\left(\frac{\alpha+1}{\delta t}-\frac{L}{2}\right)\|x^{n+1}-x^n\|^2+\|x^{n+1}-x^n\|_M^2+\langle M(x^{n+1}-x^n),-\beta_n(x^n-x^{n-1})\rangle\\
        &\geq \left(\frac{\alpha+1}{\delta t}-\frac{L}{2}\right)\|x^{n+1}-x^n\|^2+(1-\frac{\beta_n}{2})\|x^{n+1}-x^n\|_M^2-\frac{\beta_n}{2}\|x^n-x^{n-1}\|_M^2
    \end{align*}
    which leads to \eqref{lowerbound}.
\end{proof}

The following lemma lays a significant foundation for the subsequent convergence proof.
\begin{lemma}\label{lem:square:summable}
Let $\{x^n\}_n$ be a sequence generated from Algorithm \ref{alg:cap}. With the condition $\delta t<\frac{8}{77L}$, the sequence $\{\|x^{n+1}-x^n\|\}_n$ is square summable, i.e. 
\begin{align}\label{squaresum}
    \sum_{n=1}^{\infty}\|x^{n+1}-x^n\|^2<\infty.
\end{align}
\end{lemma}
\begin{proof}
 From \eqref{lowerbound}, we have derived a lower bound of $E^n(x^n)-E^n(x^{n+1})$. Now, let us turn to the corresponding upper bound. With the definition of $E^n(x)$ and $E(x)$, we have
\begin{align*}
   &\quad E^n(x^n)-E^n(x^{n+1})\\
    &= [H(x^n)+\frac{4}{11\delta t}\langle x^n-x^{n-1},x^{n-1}-x^{n-2}\rangle+F(x^n)+2\langle f(x^n)-f(x^{n-1}),x^n-x^{n-1}\rangle\\
    &\quad-\langle f(x^{n-1})-f(x^{n-2}),x^n-x^{n-2}\rangle]-[H(x^{n+1})+\frac{1}{\delta t}\|x^{n+1}-x^n\|^2-\frac{\alpha}{\delta t}\|x^{n+1}-x^n\|^2\\
    &\quad-\frac{14}{11\delta t}\langle x^{n+1}-x^n,x^n-x^{n-1}\rangle +\frac{4}{11\delta t}\langle x^{n+1}-x^{n-1},x^{n-1}-x^{n-2}\rangle+F(x^{n+1})\\
    &\quad+2\langle f(x^n)-f(x^{n-1}),x^{n+1}-x^{n-1}\rangle-\langle f(x^{n-1})-f(x^{n-2}),x^{n+1}-x^{n-2}\rangle ]\\
     &= E(x^n)-E(x^{n+1})+\frac{\alpha-1}{\delta t}\|x^{n+1}-x^n\|^2+\frac{4}{11\delta t}\langle x^n-x^{n+1},x^{n-1}-x^{n-2}\rangle\\
    &\quad +2\langle f(x^n)-f(x^{n-1}),x^n-x^{n+1}\rangle-\langle f(x^{n-1})-f(x^{n-2}),x^n-x^{n+1}\rangle] \\
    &\quad +\frac{14}{11\delta t}\langle x^{n+1}-x^n,x^n-x^{n-1}\rangle. 
\end{align*}
Utilizing the Lipschitz continuity of $f(x)$, it is clear that for any $a,b,c,d\in X$,
\begin{align*}
    \langle f(a)-f(b),c-d\rangle\leq \frac{L}{2}(\|a-b\|^2+\|c-d\|^2),
\end{align*}
which in turn presents
    \begin{align*}
    & E^n(x^n)-E^n(x^{n+1})\\
    \leq& E(x^n)-E(x^{n+1})+\frac{\alpha-1}{\delta t}\|x^{n+1}-x^n\|^2+\frac{2}{11\delta t}(\|x^n-x^{n+1}\|^2+\|x^{n-1}-x^{n-2}\|^2)\\
    &+L(\|x^n-x^{n-1}\|^2+\|x^{n}-x^{n+1}\|^2)+\frac{L}{2}(\|x^{n-1}-x^{n-2}\|^2+\|x^{n+1}-x^{n}\|^2)\\
    &+\frac{7}{11\delta t}(\|x^{n+1}-x^{n}\|^2+\|x^{n}-x^{n-1}\|^2)\\
    =&E(x^n)-E(x^{n+1})+\left(\frac{11\alpha-2}{11\delta t}+\frac{3}{2}L\right)\|x^{n+1}-x^n\|^2+\left(\frac{7}{11\delta t}+L\right)\|x^n-x^{n-1}\|^2\\
    &+\left(\frac{2}{11\delta t}+\frac{L}{2}\right)\|x^{n-1}-x^{n-2}\|^2.
\end{align*}
We then propose to reformulate the right-hand side of the inequality as the difference of uniformly perturbed energy across two stages, with the intention of employing recursive relations to substantiate the boundedness in subsequent proofs. To achieve this, we categorize the original energy and quadratic terms as
\begin{align*}
     &E^n(x^n)-E^n(x^{n+1})\\
     \leq &[E(x^n)+\left(\frac{9}{11\delta t}+\frac{3}{2}L\right)\|x^n-x^{n-1}\|^2+\left(\frac{2}{11\delta t}+\frac{L}{2}\right)\|x^{n-1}-x^{n-2}\|^2]\\
     &-[E(x^{n+1})+\left(\frac{9}{11\delta t}+\frac{3}{2}L\right)\|x^{n+1}-x^{n}\|^2+\left(\frac{2}{11\delta t}+\frac{L}{2}\right)\|x^{n}-x^{n-1}\|^2]\\
     &+\left(\frac{11\alpha+7}{11\delta t}+3L\right)\|x^{n+1}-x^n\|^2.
\end{align*}
Together with~\eqref{lowerbound}, we have the following inequality
\begin{align*}
    &\left(\frac{\alpha+1}{\delta t}-\frac{L}{2}\right)\|x^{n+1}-x^n\|^2+(1-\frac{\beta_n}{2})\|x^{n+1}-x^n\|_M^2-\frac{\beta_n}{2}\|x^n-x^{n-1}\|_M^2 \\
    \leq &[E(x^n)+\left(\frac{9}{11\delta t}+\frac{3}{2}L\right)\|x^n-x^{n-1}\|^2+\left(\frac{2}{11\delta t}+\frac{L}{2}\right)\|x^{n-1}-x^{n-2}\|^2]\\
     &-[E(x^{n+1})+\left(\frac{9}{11\delta t}+\frac{3}{2}L\right)\|x^{n+1}-x^{n}\|^2+\left(\frac{2}{11\delta t}+\frac{L}{2}\right)\|x^{n}-x^{n-1}\|^2]\\
     &+\left(\frac{11\alpha+7}{11\delta t}+3L\right)\|x^{n+1}-x^n\|^2.
\end{align*}
By \(\{\beta_n\}_n \subset [0, 1)\), we get
\begin{align*}
    &\frac{1-\beta_n}{2} \|x^{n+1} - x^n\|_M^2 +\left(\frac{4}{11\delta t}-\frac{7}{2}L\right)\|x^{n+1}-x^n\|^2 \\
    \leq &[E(x^n)+\left(\frac{9}{11\delta t}+\frac{3}{2}L\right)\|x^n-x^{n-1}\|^2+\left(\frac{2}{11\delta t}+\frac{L}{2}\right)\|x^{n-1}-x^{n-2}\|^2+\frac{1}{2}\|x^n-x^{n-1}\|_M^2]\\
     &-[E(x^{n+1})+\left(\frac{9}{11\delta t}+\frac{3}{2}L\right)\|x^{n+1}-x^{n}\|^2+\left(\frac{2}{11\delta t}+\frac{L}{2}\right)\|x^{n}-x^{n-1}\|^2+\frac{1}{2}\|x^{n+1}-x^{n}\|_M^2].
\end{align*}
Now denote 
\begin{align}\label{axyz}
    A(x,y,z)=E(x)+\left(\frac{9}{11\delta t}+\frac{3L}{2}\right)\|x-y\|^2+\left(\frac{2}{11\delta t}+\frac{L}{2}\right)\|y-z\|^2+\frac{1}{2}\|x-y\|_M^2,
\end{align}
and we have 
\begin{align}\label{squarebound}
    \frac{1-\beta_n}{2} \|x^{n+1} - x^n\|_M^2 +\left(\frac{4}{11\delta t}-\frac{7}{2}L\right)\|x^{n+1}-x^n\|^2 \leq A(x^n,x^{n-1},x^{n-2})-A(x^{n+1},x^n,x^{n-1}).
\end{align}
Summing the inequality above from 2 to $N$, we have
\begin{align}
    \sum_{n=2}^N (M_1\|x^{n+1} - x^n\|_M^2+M_2\|x^{n+1}-x^n\|^2)\leq A(x^2,x^{1},x^{0})-A(x^{N+1},x^N,x^{N-1})
\end{align}
where 
\begin{align*}
    M_1=\frac{1-\beta_n}{2}>0,\quad M_2=\frac{4}{11\delta t}-\frac{7}{2}L>0
\end{align*}
under the conditions $\delta t<\frac{8}{77L}$ and $\{\beta_n\}_n \subset [0, 1)$. Taking $N\rightarrow\infty$, we get
\begin{align}
    \sum_{n=1}^{\infty}\|x^{n+1}-x^{n}\|^2<\infty.
\end{align}
This completes the proof.
\end{proof}

Now we turn to the global convergence. We first demonstrate that \(\text{dist}(\mathbf{0}, \partial A(x^{n+1}, x^n, x^{n-1}))\) is bounded above by a sum of consecutive differences \(\|x^n - x^{n-1}\|\). Then, since lemma \ref{lem:square:summable} shows that \(A(x^{n+1}, x^n, x^{n-1}) - A(x^{n+2}, x^{n+1}, x^n)\) is bounded below by \(\|x^{n+2} - x^{n+1}\|^2\), we use the KL property to connect the estimates above. This yields \(\sum_{n=0}^{\infty} \|x^{n+1} - x^n\| < \infty\), implying that \(\{x^n\}_n\) is a Cauchy sequence and hence converges.

\begin{theorem}\label{thm:globalconverg}
  Assuming that $f(x)$ is Lipschitz continuous with parameter $L$ and $A(x,y,z)$ is a KL function generated by \eqref{axyz} from Algorithm \ref{alg:cap}, then the following statements hold:
\begin{itemize}
    \item[\rm (i)] $\lim\limits_{n\rightarrow\infty} {\dist} (\emph{\tbfzero},\partial A(x^{n+1},x^n,x^{n-1}))=0$.
    \item[\rm (ii)] $\{A(x^{n+1},x^n,x^{n-1})\}_n$ is non-increasing with a limit and $\{x^n\}_n$ is bounded. Specifically, there exists a constant $\zeta$, such that $\lim\limits_{n\rightarrow\infty}A(x^{n+1},x^n,x^{n-1})=\zeta$.
    \item[\rm (iii)] $\{x^n\}_n$ converges to a critical point of $E$ and $\sum_{n=0}^{\infty}\|x^{n+1}-x^n\|<\infty$.
\end{itemize}

\end{theorem}

\begin{proof}
 (i) With direct conduction, we obtain
\begin{align*}
\partial A(x,y,z) \big|_{(x,y,z) = (x^{n+1}, x^n, x^{n-1})} =
    \begin{pmatrix}
\partial H(x^{n+1}) + f(x^{n+1}) + \left(C_1I+M\right)(x^{n+1} - x^n) \\
-\left(C_1I+M\right)(x^{n+1} - x^n) + C_2(x^n - x^{n-1}) \\
-C_2(x^n - x^{n-1})
\end{pmatrix}
\end{align*}
where $C_1=\frac{18}{11\delta t} + 3L$, $C_2=\frac{4}{11\delta t} + L$. According to the first-order optimality condition, we obtain
\begin{align*}
    0\in & \partial H(x^{n+1})  +\frac{12}{11\delta t} \left[\frac{11}{6}x^{n+1} - 3x^{n} + \frac{3}{2}x^{n-1} - \frac{1}{3}x^{n-2}\right] + M(x^{n+1} - y^{n})\\
    & + [3f(x^{n}) - 3f(x^{n-1}) + f(x^{n-2})].
\end{align*}
Then $\text{dist}(\textbf{0},\partial  A(x^{n+1}, x^n, x^{n-1}))$ can be estimated as
\begin{align*}
    &\text{dist}(\textbf{0},\partial A(x^{n+1}, x^n, x^{n-1}))\\
    \leq& \|f(x^{n+1}) - [3f(x^{n}) - 3f(x^{n-1}) + f(x^{n-2})]+ (C_1I+M)(x^{n+1} - x^n)\\
    &- \frac{12}{11\delta t} \left[\frac{11}{6}x^{n+1} - 3x^{n} + \frac{3}{2}x^{n-1}- \frac{1}{3}x^{n-2}\right]- M(x^{n+1} - y^{n}) \|\\
    &+ \|(C_1 I+M)(x^{n+1} - x^n) - C_2(x^n - x^{n-1})\| + \|C_2(x^n - x^{n-1})\| .
\end{align*}
Denote $\lambda_M$ as the largest eigenvalue of $M$. We further obtain
\begin{align*}
  &\text{dist}(\textbf{0},\partial A(x^{n+1}, x^n, x^{n-1}))\\
\leq& \left(L+2C_1+3\lambda_M+\frac{2}{\delta t}\right)\|x^{n+1} - x^{n}\|+ \left(2L+2C_2+\frac{14}{11\delta t}+\lambda_M\beta_{n}\right)\|x^{n} - x^{n-1}\|\\
&+\left(L+\frac{4}{11\delta t}\right)\|x^{n-1}-x^{n-2}\|\\
\leq & \tilde{C} (\|x^{n+1} - x^{n}\| + \|x^{n} - x^{n-1}\|+\|x^{n-1} - x^{n-2}\|)
\end{align*}
where $\tilde{C}  = \max \left\{L+2C_1+3\lambda_M+\frac{2}{\delta t},  2L+2C_2+\frac{14}{11\delta t}+\lambda_M\beta_{n}, L+\frac{4}{11\delta t}\right\}$.

Since $\|x^{n+1}-x^n\| \to 0$ as $n \to \infty$, we have
\begin{align*}
\text{dist}(\textbf{0},\partial A(x^{n+1}, x^n, x^{n-1}))\rightarrow0 \quad (n \to \infty).
\end{align*}

(ii) With \eqref{squarebound}, $\{A(x^{n+1}, x^n, x^{n-1})\}_n$ is a non-negative and non-increasing sequence, so the sequence is bounded and there is a limit $\zeta$. With the level-boundedness of $E(x)$, $A(x, y, z)$ would be level-bounded as well. Since $ A(x^{n+1}, x^n, x^{n-1})$ is bounded and level-bounded on $(x^{n+1}, x^n, x^{n-1})$, the sequence $\{x^n\}_n$ is bounded.  

(iii) Assuming that $\psi(\cdot)$ is a concave function, we have
\begin{equation}\label{ineq:KL_of_psi}
    \begin{aligned}
&[\psi(A(x^{n+1}, x^n, x^{n-1}) - \zeta) - \psi(A(x^{n+2}, x^{n+1}, x^n) - \zeta)] \cdot\text{dist}(\textbf{0},\partial A(x^{n+1}, x^n, x^{n-1}))\\
\geq &\psi'(A(x^{n+1}, x^n, x^{n-1}) - \zeta)[A(x^{n+1}, x^n, x^{n-1}) - A(x^{n+2}, x^{n+1}, x^n)]\\
&\cdot\text{dist}(\textbf{0},\partial A(x^{n+1}, x^n, x^{n-1}))\\
= &\underbrace{\psi'(A(x^{n+1}, x^n, x^{n-1}) - \zeta)\cdot\text{dist}(\textbf{0},\partial A(x^{n+1}, x^n, x^{n-1}))}_{\geq 1}\\
&\cdot [A(x^{n+1}, x^n, x^{n-1}) - A(x^{n+2}, x^{n+1}, x^n)]\\
 \geq &A(x^{n+1}, x^n, x^{n-1}) - A(x^{n+2}, x^{n+1}, x^n)\\
\geq &M_2 \|x^{n+2} - x^{n+1}\|^2
\end{aligned}
\end{equation}
where the second inequality holds with the KL-property of  $A(x, y, z)$.

Now denote 
\begin{equation}
    \phi(x^{n-1}, x^n, x^{n+1}, x^{n+2}, \zeta)= \psi(A(x^{n+1}, x^n, x^{n-1}) - \zeta) - \psi(A(x^{n+2}, x^{n+1}, x^n) - \zeta).
\end{equation}
Then we have
\begin{align*}
\|x^{n+2} - x^{n+1}\|^2 \leq \tilde{K} (\|x^{n+2} - x^{n+1}\| + \|x^{n+1} - x^n\| + \|x^n - x^{n-1}\|)\phi(x^{n-1}, x^n, x^{n+1}, x^{n+2}, \zeta)
\end{align*}
where $\tilde{K} = \frac{1}{M_2}\tilde{C}$. Furthermore, by the following deduction
\begin{align*}
a^2 \leq cd \quad \Rightarrow \quad a \leq c + \frac{d}{4}
\end{align*}
we have
\begin{align*}
\|x^{n+2} - x^{n+1}\| \leq \tilde{K} \phi(x^{n-1}, x^n, x^{n+1}, x^{n+2}, \zeta) + \frac{1}{4} (\|x^{n+2} - x^{n+1}\| + \|x^{n+1} - x^n\| + \|x^n - x^{n-1}\|)
\end{align*}
which is equivalent to
\begin{align*}
\frac{1}{4} \|x^{n+2} - x^{n+1}\| \leq \tilde{K} \phi(x^{n-1}, x^n, x^{n+1}, x^{n+2}, \zeta) + \frac{1}{4} (\|x^{n+1} - x^n\| + \|x^n - x^{n-1}\| - 2\|x^{n+2} - x^{n+1}\|).
\end{align*}
Denoting $\nu_n=\|x^{n+1} - x^{n}\|$, we conclude that
\begin{align}
\nu_{n+1} \leq 4\tilde{K} \phi(x^{n-1}, x^n, x^{n+1}, x^{n+2}, \zeta) + (\nu_n - \nu_{n+1})+ (\nu_{n-1} - \nu_{n+1})
 \end{align}
 i.e. 
 \begin{align}\label{ineq:splitform}
\nu_{n} \leq 4\tilde{K} \phi(x^{n-2}, x^{n-1}, x^n, x^{n+1}, \zeta) + (\nu_{n-1} - \nu_{n})+ (\nu_{n-2} - \nu_{n}).
 \end{align}
Since $\nu_n\to 0$, the summations of the two series in the inequality are respectively:
\begin{align*}
    &\sum_{n=T}^{\infty} (\nu_{n-1} - \nu_{n})=\nu_{T-1},\\
    &\sum_{n=T}^{\infty} (\nu_{n-2} - \nu_{n})=\nu_{T-2}+\nu_{T-1}.
\end{align*}
Summing the inequality from $T$ to $\infty$, we have 
\begin{align*}
\sum_{n=T}^{\infty} \nu_{n} \leq 4\tilde{K} \psi(A(x^T, x^{T-1}, x^{T-2}) - \zeta) + \nu_{T-1}+ (\nu_{T-2} + \nu_{T-1}).
\end{align*}
It is equivalent to
\begin{align*}
\sum_{n=T}^{\infty} \|x^{n+1} - x^n\| \leq 4\tilde{K} \psi(A(x^T, x^{T-1}, x^{T-2}) - \zeta) + 2\|x^T - x^{T-1}\| + \|x^{T-1} - x^{T-2}\| < +\infty
\end{align*}
This implies that\(\{x^n\}_n\) is a Cauchy sequence and hence converges. 
\end{proof}

In what follows, we will demonstrate the analysis of local convergence.

\subsection{Local convergence}
We need the KL exponent of the $A(x,y,z)$ for the local convergence rate. The following theorem and its proof are standard (see \cite[Theorem 2]{Attouch2009} or \cite[Lemma 1]{Artacho2018}). 

\begin{theorem}[Local convergence rate]\label{them:local:rate}Under the same assumption as Theorem \ref{thm:globalconverg}, consider a sequence $\{x^n\}_n$ produced by Algorithm \ref{alg:cap} that converges to $x^*$. Assume that $A(x,y,z)$ is a $KL$ function with $\phi$ in the $KL$ inequality given as $\phi(s) = cs^{1-\theta}$, where $\theta\in[0,1)$ and $c>0$. The following statements are valid:
\begin{itemize}
 \item[\emph{(i)}] For $\theta = 0$, there exists a positive $n_0$ such that $x^n$ remains constant for $n > n_0$.
\item[\emph{(ii)}] For $\theta\in (0, \frac12]$, there are positive constants $c_1$, $n_1$, and $\eta\in (0,1)$ such that $\|x^n - x^*\| < c_1\eta^n$ for $n > n_1$.
\item[\emph{(iii)}] For $\theta\in (\frac{1}{2},1)$, there exist positive constants $c_2$ and $n_2$ such that $\|x^n-x^*\| < c_2 n^{-\frac{1-\theta}{2\theta-1}} $ for $n>n_2$.
\end{itemize}
\end{theorem}
\begin{proof}
(i) If $\theta = 0$, we claim that there must exist $n_0 > 0$ such that
\begin{align*}
    A(x^{n_0+1},x^{n_0},x^{n_0-1}) = \zeta.
\end{align*}
 In fact, suppose to the contrary that  $A(x^{n+1},x^{n},x^{n-1}) > \zeta$ for all $n>0$. Since $\lim\limits_{n\to \infty}x^n = x^*$ and the sequence $\{A(x^{n+1},x^n,x^{n-1})\}_n$ is monotone decreasing and convergent to $\zeta$ by Theorem \ref{thm:globalconverg}(ii), we choose the concave function $\psi(s) = cs$ and the KL inequality \eqref{eq:kl:def} that for all sufficiently large $n$,
 \begin{align*}
\text{dist}(\textbf{0},\partial A(x^{n+1}, x^n, x^{n-1})) \geq c^{-1}
\end{align*}
which contradicts \rm{Theorem \ref{thm:globalconverg}\rm{(i)}}. Thus, there exist $n_0 > 0$ so that $A(x^{n_0+1},x^{n_0},x^{n_0-1}) = \zeta$.
 
 Since the sequence $\{A(x^{n+1},x^n,x^{n-1})\}_n$ is monotone decreasing and convergent to $\zeta$, it must hold that $A(x^{n_0+\Bar{n}},x^{n_0+\Bar{n}-1}) = \zeta$ for any $\Bar{n}>0$. Thus, we can conclude from \eqref{ineq:KL_of_psi} that $x^{n_0} = x^{n_0+\Bar{n}}$. This proves that if $\theta = 0$, there exists $n_0>0$ so that $x^n$ is constant for $n>n_0$.

(ii) When $\theta \in (0,1)$, we only need to consider the case when $A(x^n, x^{n-1},x^{n-2})>\zeta$ for all $n>0$ based on the proof above.

Define $\Delta_n =A(x^{n},x^{n-1},x^{n-2}) - \zeta$ and $S_n = \sum_{i=n}^{\infty}\|x^{i+1} - x^{i}\|$, where $S_n$ is well-defined due to the Theorem \ref{thm:globalconverg}(iii). Then, from (\ref{ineq:splitform}), we have for any $n>N$ that 
\begin{align*}
    S_n = &\sum_{i=n}^{\infty}\|x^{i+1} - x^i\|\\
    \leq &\sum_{i=n}^{\infty}\Big[4\tilde{K} \phi(x^{n-2}, x^{n-1}, x^n, x^{n+1}, \zeta) + [\|x^n - x^{n-1}\| - \|x^{n+1} - x^n\|]\\
&\quad \quad + [\|x^{n-1} - x^{n-2}\| - \|x^{n+1} - x^n\|]\Big]\\
    \leq&4\tilde{K}\psi(A(x^{n},x^{n-1}, x^{n-2}) - \zeta)  + 2\|x^n - x^{n-1}\| + \|x^{n-1} - x^{n-2}\|\\
    =&4\tilde{K}\psi(\Delta_n)  + 2(S_{n-1}-S_{n})+(S_{n-2}-S_{n-1})\\
    \leq&4\tilde{K}\psi(\Delta_n)  + 2(S_{n-2}-S_{n+1}).
\end{align*}

Set $ \psi(s) = cs^{1-\theta}$. For all sufficiently large $n$, we have
\begin{equation*}
   c(1-\theta)\Delta_n^{-\theta}\cdot\text{dist}(\textbf{0},\partial A(x^{n}, x^{n-1}, x^{n-2}))\geq 1.
\end{equation*} 
We can transfer the estimate of $\text{dist}(\textbf{0},\partial A(x^{n},x^{n-1},x^{n-2}))$ to a new formulation
\begin{equation*}
   \text{dist}(\textbf{0},\partial A(x^{n},x^{n-1},x^{n-2}))\leq  \tilde{C} (\|x^{n+1} - x^{n}\| + \|x^{n} - x^{n-1}\|+\|x^{n-1} - x^{n-2}\|)= \tilde{C}(S_{n-2}-S_{n+1}).
\end{equation*}
Due to the above two inequalities, we have
\begin{equation*}
    \Delta_n^{\theta}\leq c(1-\theta) \tilde{C}(S_{n-2}-S_{n+1}).
\end{equation*}
Combining it with $S_n\leq4\tilde{K}\psi(\Delta_n)  + 2(S_{n-2}-S_{n+1})$, we have 
\begin{align}\label{eq:ineq_Sn}
    S_n\leq& 4c\tilde{K}(\Delta_n^{\theta})^{\frac{1-\theta}{\theta}} + 2(S_{n-2}-S_{n+1})\leq C_1(S_{n-2}-S_{n+1})^{\frac{1-\theta}{\theta}} + 2(S_{n-2}-S_{n+1})
\end{align}
where $C_1 = 4c\tilde{K}[c(1-\theta)\tilde{C}]^{\frac{1-\theta}{\theta}}$. 

Suppose the first that $\theta \in (0,\frac12]$. We have $\frac{1-\theta}{\theta} \geq 1$. Since $\|x^{n+1} - x^n\| \to 0$ from Lemma \ref{lem:square:summable}, it leads to $S_{n-2} - S_{n+1} \to 0$. From these and \eqref{eq:ineq_Sn}, we can conclude that there exists $n_1>0$ such that for all $n\geq n_1$, we have
\begin{equation*}
    S_{n+1}\leq S_n\leq (C_1+2)(S_{n-2}-S_{n+1})
\end{equation*}
which implies that 
\begin{equation*}
    S_{n+1} \leq \frac{C_1+2}{C_1+3}S_{n-2}.
\end{equation*} 
Hence, for sufficiently large $n>n_1$,
\begin{equation*}
    \|x^{n+1} - x^*\|\leq \sum_{i = n+1}^{\infty}\|x^{i+1} - x^i\| = S_{n+1}\leq S_{n_1-2}\eta^{n-n_1}, \quad  \eta := \sqrt{\frac{C_1 + 2}{C_1 + 3}}.
\end{equation*}
(iii) Finally, we consider the case that $\theta \in (\frac12,1)$, which imply $\frac{1-\theta}{\theta}\leq 1$. Combining this with the fact that $S_{n-2} - S_n \to 0$, we see that there exists $n_2 > 0$ such that for all $n\geq n_2$, we have
\begin{align*}
    S_{n+1}\leq S_n\leq& 4c\tilde{K}(\Delta_n^{\theta})^{\frac{1-\theta}{\theta}} + 2(S_{n-2}-S_{n+1})
    =(C_1 + 2)(S_{n-2}-S_{n+1})^{\frac{1-\theta}{\theta}}.
\end{align*}
Raising both sides of the above inequality to the power of $\frac{\theta}{1-\theta}$, we can observe a new inequality that 
\begin{align*}
 S_{n+1}^{\frac{\theta}{1-\theta}} \leq C_2(S_{n-2}-S_{n+1}) 
\end{align*}
for $n \geq n_2$, where $C_2 = (C_1 + 2)^{\frac{\theta}{1-\theta}}$. Let us define the sequence $\Omega_n = S_{3n}$. For any $n \geq \left\lceil \frac{n_2}{2} \right\rceil$, with nearly the same arguments as in \cite[Page 15]{Attouch2009}, there exists a constant $C_3>0$ such that for sufficient large $n$ 
\begin{equation*}
    \Omega^{\frac{\theta}{1-\theta}}_{n} \leq C_2(\Omega_{n-1}- \Omega_{n}) \Rightarrow
    \Omega_{n} \leq C_3n^{-\frac{1-\theta}{2\theta -1}}.
\end{equation*}
Consequently, for any positive integer $n$,
\begin{align*}
    \|x^{n+1} - x^*\|
    \leq S_{n+1}\left\{
    \begin{array}{ll}
         = \Omega_{\frac{n}{3}}\leq 3^{\rho}C_3n^{-\rho} &\text{if $3|n$}\\
         \leq \Omega_{\frac{n-1}{3}}\leq 3^{\rho}C_3(n-1)^{-\rho}\leq 6^{\rho}C_3n^{-\rho}   &\text{if $3|n-1$}\\
         \leq \Omega_{\frac{n-2}{3}}\leq 3^{\rho}C_3(n-2)^{-\rho}\leq 6^{\rho}C_3n^{-\rho}   &\text{if $3|n-2$}
    \end{array}
    \right.
\end{align*}
where $ \rho := \frac{1-\theta}{2\theta-1}$. The proof is completed. 
\end{proof}

At this point, we conclude the convergence analysis of Algorithm \ref{alg:cap}, encompassing both global and local convergence properties. In the next section, we will present numerical experiments to demonstrate the efficiency of our algorithm.

\section{Numerical experiments} \label{sec:pre}
In this section, we perform two numerical experiments to demonstrate the efficiency of our algorithm in solving nonconvex optimization problems, which are respectively.

\textbf{(1) Least squares problems with (modified) SCAD regularizer}: conducted on a computer with Intel(R) Core(TM) i5-10210U CPU @ 2.11 GHz.

\textbf{(2) Graphic Ginzburg-Landau model}: executed on a workstation with Intel(R) Xeon(R) CPU E5-2699A v4 (2.40GHz) and a GPU of NVIDIA GeForce RTX 2080 Ti (22GB).

Implementation details, parameter settings, and comparative results will be discussed below. For the KL properties of the two problems, please refer to Remark \ref{KL:two:models}.

\subsection{Least squares problems with (modified) SCAD regularizer}
 We consider the smoothly clipped absolute deviation (SCAD) regularization, whose DC decomposition can be expressed as (see \cite[Section 6.1]{APX} or \cite{Wen2018})
 \begin{equation*}
    P(x) = \mu\sum_{i=1}^{k}\int_0^{|x_i|}\min\left\{1,\frac{[\theta\mu-x]_+}{(\theta-1)\mu}\right\}dx=\mu\|x\|_1 - \underbrace{\mu\sum_{i=1}^{k}\int_0^{|x_i|}\frac{[\min\{\theta\mu,x\}-\mu]_+}{(\theta-1)\mu}dx}_{\tilde P(x)}
\end{equation*}
 where $\theta>2$ is a constant, $\mu>0$ serves as the regularization parameter, $[x]_{+}=\max\{0,x\}$ and $\tilde P(x) =  \mu \|x\|_1 -P(x) = \sum_{i=1}^k \tilde p_i(u_i)$. One can easily verify that $\tilde P(x)$ is continuously differentiable with the gradient
 \begin{equation*}
    \tilde  p_i(x_i) = \left\{
        \begin{array}{ll}
         0& \text{if}  \ \  \ |x_i|\leq\mu \\
         \frac{(|x_i| - \mu)^2}{2(\theta - 1)\mu}& \text{if} \ \  \ \mu<|x_i|<\theta \\
         \mu|x_i| - \frac{\mu^2(\theta+1)}{2}&\text{if}\ \ \ |x_i|\geq\theta\mu 
    \end{array},
        \right. \ \ \nabla_i \tilde P_i(x_i)= \text{sign}(x_i) \dfrac{[\min\{\theta \mu,|x_i|\} -\mu]_{+}} {(\theta-1)}.
\end{equation*}
 
 Applying SCAD regularization in the least squares problem, we obtain the following optimization formulation:
\begin{equation}\label{scad}
    \min_{x\in \mathbb{R}^k}E(x) = \frac12\|Ax-b\|^2 + P(x),
\end{equation}
where $A\in \mathbb{R}^{m\times k}$, $b\in \mathbb{R}^{m}$.
 
Also, for our smooth case, we consider a variant of SCAD regularization. To achieve that, we substitute the $l_1$ norm with the Huber function, which is denoted as 
\begin{align*}
    \mathcal{H}_{\gamma}(x):= \sum_{i}\mathcal{H}(x_i,\gamma)
\end{align*}
where $\mathcal{H}(u_i,\gamma)$ has the form that
\begin{equation*}
    \mathcal{H}(x_i,\gamma) = \left\{
    \begin{array}{ll}
         \frac{|x_i|^2}{2\gamma}& \text{if}  \ \  \ |x_i|\leq\gamma \\
         |x_i|-\frac{\gamma}{2}& \text{if} \ \  \ |x_i|>\gamma
    \end{array}.
    \right. 
\end{equation*}
It is worth noting that $\gamma$ is a shape parameter characterizing the level of robustness, and we set $\gamma=0.5\mu$ in the following experiment. Incorporating the modified SCAD regularization into the least squares problem, we obtain
\begin{equation}\label{huberscad}
    \min_{x\in \mathbb{R}^k}E(x) = \frac12\|Ax-b\|^2 + \mu \mathcal{H}_{\gamma}(x) - \tilde P(x),
\end{equation}
where $A\in \mathbb{R}^{m\times k}$, $b\in \mathbb{R}^{m}$. Note that the modified SCAD regularization $P_{M}(x): =  \mu\mathcal{H}_{\gamma}(x) - \tilde P(x)$ admits a separable structure:
\begin{equation*}
    P_{M}(x) := \sum_{i=1}^k p_{M,i}(x_i),
\end{equation*} 
where each component function $p_{M,i}(x_i)$ is given by the following piecewise definition:
        \begin{equation}\label{eq:huber-scad}
        \frac{p_{M,i}(x_i)}{\mu} = \left\{
        \begin{array}{ll}
       |x_i|^2/(2\gamma)& \text{if}  \ \  \ |x_i|\leq\gamma \\
         |x_i|-{\gamma}/{2}& \text{if}  \ \  \ \gamma < |x_i|\leq\mu \\
         |x_i|-{\gamma}/{2} -(|x_i| - \mu)^2/(2(\theta - 1)\mu)& \text{if} \ \  \ \mu<|x_i|<\theta\mu \\
         (\mu(\theta+1) - \gamma)/{2} &\text{if}\ \ \ |x_i|\geq\theta\mu 
    \end{array}.
        \right.
    \end{equation}
To formulate the least square problem, we first generate an $m \times k$ random matrix $A$ with normalized columns ($\|A_j\|_2 = 1$ for all $j = 1,...,k$). We then construct a sparse vector $y \in \mathbb{R}^k$ by uniformly randomly selecting its support $T \subset \{1,...,k\}$ with size $|T| = s$. In other words, $s$ is designed as the number of non-zero elements in $y$, which characterizes its sparsity level. The vector $b \in \mathbb{R}^m$ is generated according to:
\begin{equation*}
    b = Ay +  0.01 \xi,\quad \xi \sim \mathcal{N}(0,I_m),
\end{equation*}
where $\xi$ consists of i.i.d. standard Gaussian entries. All algorithms are initialized at the origin and terminate when the relative step difference satisfies:
\begin{equation*}
    \frac{\|x^{n} - x^{n-1}\|}{\max\{1, \|x^{n}\|\}} < 10^{-12}.
\end{equation*}

  We focus on problem sizes characterized by the tuple \( (m, k, s) = (720i, 2560i, 80i) \), where  \( i \) is an integer ranging from 1 to 10. This setup aligns with the high-dimensional setting in sparse statistics where the data dimension exceeds the number of data points, and the solution is sufficiently sparse \cite{hastie2015statistical}. 

In our experiments, we compute normal SCAD \eqref{scad} and Huber-SCAD \eqref{huberscad} with different algorithms for comparative analysis. In the evaluation of various algorithms, the iteration count (denoted as \( \text{iter} \)) and the CPU time are meticulously recorded. Concurrently, we record the number of nonzero elements in the output vector to ascertain whether the solution adheres to the prescribed sparsity constraints.

The computational results are summarized in Table~\ref{table1} and Table~\ref{table2}, corresponding to problems~\eqref{scad} and~\eqref{huberscad} with parameters $\mu = 0.033$ and $\theta = 10$. For each problem size, five random instances are generated, and the values reported in each row of Table~\ref{table1} and Table~\ref{table2} represent the average performance across all instances. We compare our proposed algorithm, 3BapDCA$_{e}$, with four established algorithms: 3BapDCA, BapDCA, DC algorithm(DCA), and BDCA - Backtracking algorithm as described in~\cite{Artacho2018}. The implementation details of these algorithms are discussed below.

\begin{itemize}
    \item \textbf{3BapDCA$_{e}$}: This algorithm is represented by Algorithm \ref{alg:cap} where $H(x) = \mu \mathcal{H}_{\gamma}(x) + \frac12\|Ax-b\|^2$ and $F(x) =  - \tilde P(x)$. The Lipschitz constant $L$ of $\nabla F$ is $1/(\theta-1)$ \cite[Example 4.3]{Wen2018} with $L=1/9$. The parameter values for this algorithm are chosen as follows:  $\alpha =2$, $\delta t < 8/(77L)$ (e.g., $\delta t= 72/77-10^{-15}$),  $\gamma = 0.2$ and $M = \lambda_{A^TA} I- A^TA$, where $\lambda_{A^TA}$ represents the largest eigenvalue of the matrix $A^TA$. The extrapolation parameters are chosen as in \eqref{extrapara} and the restart strategy is described as in section 2. 
    \item \textbf{3BapDCA}: This special version of the algorithm {3BapDCA$_{e}$}   does not include the extrapolation technique. The parameters of this algorithm remain consistent with those of the algorithm {3BapDCA$_{e}$}.
    \item \textbf{BapDCA}: This algorithm is represented by the preconditioned second-order framework \cite{shen2024secondorder}, with the energy form
    \begin{equation*}
    \begin{aligned}
           E^n(x)&=H^n(x)-F^n(x),\\
        H^n(x)&=H(x)+\frac{1}{\delta t}\|x-x^n\|^2,\\
        F^n(x)&=\frac{1}{3\delta t}\|x-x^{n-1}\|^2-F(x)-\langle f(x^n)-f(x^{n-1}),x-x^{n-1}\rangle
    \end{aligned}
    \end{equation*}
    where $H(x) = \mu \mathcal{H}_{\gamma}(x) + \frac12\|Ax-b\|^2$ and $F(x) =  - \tilde P(x)$. The Lipschitz constant $L$ of $\nabla F$ is $1/(\theta-1)$ \cite[Example 4.3]{Wen2018} with $L=1/9$. The step size satisfies $\delta t < 2/(3L)$ (e.g., $\delta t= 6-10^{-15}$). We set $\gamma = 0.2$ and $M = \lambda_{A^TA} I- A^TA$ as well. 
    \item \textbf{BDCA}: This algorithm is based on a combination of DCA together with a line search technique, as detailed in \cite[Algorithm 2]{Artacho2018}. To simplify computations, we employ distinct convex splitting to sidestep equation solving. The splitting can be defined as follows:
    \begin{equation*}
        E(x) =\mu \mathcal{H}_{\gamma}(x) + \frac{\lambda_{A^TA}}{2}\|x\|^2 - \left(\frac{\lambda_{A^TA}}{2}\|x\|^2 +  \tilde P(x)  - \frac12\|Ax-b\|^2\right)
    \end{equation*}
    where $\lambda_{A^TA}$ represents the maximum eigenvalue of the matrix $A^TA$. The Lipschitz constant $L$ of $\nabla F$ is $1/(\theta - 1)$ with $L = 1/9$. The parameter values of the line search part for this algorithm are chosen as follows: $\delta t < 2/(3L)$ (e.g., $\delta t = 6 \times 10^{-7}$), $\lambda_{\max} = 5$, $\lambda = 0.618\lambda_{\max}$, $\alpha = 0.2$, $\beta = 0.8$.
    \item \textbf{DCA}: This is the classical DC algorithm, which is a special version of the algorithm BDCA without the line search \cite[Algorithm 1]{Artacho2018}. 
\end{itemize}
From Table \ref{table1} and Table \ref{table2}, it is evident that our proposed 3BapDCA$_{e}$ algorithm significantly outperforms other algorithms in terms of CPU time. Although DCALS has fewer iterations in Table \ref{table1}, it takes a much longer time. Moreover, after a slight smoothing being introduced to the objective functional, DCALS costs more iterations than 3BapDCA$_{e}$ in Table \ref{table2}, which illustrates the stability of our 3BapDCA$_{e}$ in different problems. It can be seen that the preconditioning can alleviate the restrictions on step sizes and bring out more flexibility. Furthermore, with  $M = \lambda_{A^TA} I- A^TA$, one can choose $F(x) = -\tilde P(x)$, $\delta t < 2/(3L)=6$ with $L$ being the Lipschitz constant of $\tilde P$. However, for computing explicit resolvents, one has to choose $F(x) = \frac{1}{2}\|Ax-b\|^2 -\tilde P(x)$ with $\delta t < 2/(3\lambda_{A^TA}) < 1/12$ since $\lambda_{A^TA}>8$, where the step size is much smaller.

\begin{sidewaystable}[htbp!]
\centering
\caption{Solving (\ref{scad}) on random instances\\ (\textcircled{1}DCA \textcircled{2}BDCA \textcircled{3}BapDCA \textcircled{4}3BapDCA \textcircled{5}3BapDCA$_\text{e}$)}\label{table1}
\begin{tabular}{ccc|rrrrr|rrrrr|c}
\hline
\multicolumn{3}{c|}{Size} & \multicolumn{5}{c|}{iter} & \multicolumn{5}{c|}{CPU time (s)} & \multicolumn{1}{c}{Sparsity}\\ 
\hline
$m$ & $k$ & $s$ & \textcircled{1} & \textcircled{2} & \textcircled{3} & \textcircled{4} & \textcircled{5} & \textcircled{1} & \textcircled{2} & \textcircled{3} & \textcircled{4} & \textcircled{5} &$s'$ \\ \hline
720 & 2560 & 80 & 557 & \textbf{131} & 369 & 636 & 173 & 0.8 & 0.8 & 0.8 & 1.4 & \textbf{0.4} & 79 \\
1440 & 5120 & 160 & 581 & \textbf{141} & 384 & 656 & 177 & 3.3 & 3.7 & 3.3 & 5.6 & \textbf{1.5} & 160 \\
2160 & 7680 & 240 & 585 & \textbf{145} & 388 & 663 & 178 & 8.7 & 10.4 & 8.2 & 13.7 & \textbf{3.7} & 243 \\
2880 & 10240 & 320 & 588 & \textbf{147} & 391 & 667 & 179 & 20.0 & 24.7 & 19.2 & 32.6 & \textbf{8.7} & 322 \\
3600 & 12800 & 400 & 593 & \textbf{145} & 391 & 671 & 180 & 31.3 & 37.2 & 32.0 & 56.6 & \textbf{13.3} & 398 \\
4320 & 15360 & 480 & 596 & \textbf{146} & 395 & 674 & 179 & 47.7 & 59.4 & 45.4 & 75.1 & \textbf{19.4} & 480 \\
5040 & 17920 & 560 & 590 & \textbf{146} & 391 & 667 & 179 & 54.2 & 64.5 & 51.9 & 88.9 & \textbf{23.7} & 564 \\
5760 & 20480 & 640 & 590 & \textbf{144} & 391 & 667 & 178 & 79.0 & 92.9 & 77.1 & 127.8 & \textbf{35.1} & 642 \\
6480 & 23040 & 720 & 595 & \textbf{144} & 394 & 672 & 180 & 84.8 & 105.0 & 82.8 & 144.0 & \textbf{37.0} & 725 \\
7200 & 25600 & 800 & 586 & \textbf{145} & 388 & 663 & 178 & 99.7 & 118.6 & 96.8 & 168.0 & \textbf{44.5} & 798 \\ \hline
\end{tabular}

\end{sidewaystable}

\begin{sidewaystable}[htbp!]
\centering
\caption{Solving (\ref{huberscad}) on random instances\\ (\textcircled{1}DCA \textcircled{2}BDCA \textcircled{3}BapDCA \textcircled{4}3BapDCA \textcircled{5}3BapDCA$_\text{e}$)}\label{table2}
\begin{tabular}{ccc|rrrrr|rrrrr|c}
\hline
\multicolumn{3}{c|}{Size} & \multicolumn{5}{c|}{iter} & \multicolumn{5}{c|}{CPU time (s)} & \multicolumn{1}{c}{Sparsity}\\ 
\hline
$m$ & $k$ & $s$ & \textcircled{1} & \textcircled{2} & \textcircled{3} & \textcircled{4} & \textcircled{5} & \textcircled{1} & \textcircled{2} & \textcircled{3} & \textcircled{4} & \textcircled{5} &$s'$ \\ \hline
720 & 2560 & 80 & 580 & 178 & 400 & 657 & \textbf{176} & 0.7 & 2.0 & 0.7 & 1.2 & \textbf{0.3} & 81 \\
1440 & 5120 & 160 & 575 & 179 & 396 & 650 & \textbf{176} & 3.2 & 5.2 & 3.3 & 5.3 & \textbf{1.5} & 158 \\
2160 & 7680 & 240 & 598 & 200 & 412 & 676 & \textbf{179} & 8.9 & 17.2 & 8.9 & 14.5 & \textbf{3.8} & 243 \\
2880 & 10240 & 320 & 580 & 209 & 400 & 657 & \textbf{176} & 15.2 & 32.0 & 15.2 & 24.6 & \textbf{6.8} & 322 \\
3600 & 12800 & 400 & 594 & 190 & 409 & 672 & \textbf{179} & 34.7 & 65.2 & 34.3 & 56.2 & \textbf{15.3} & 401 \\
4320 & 15360 & 480 & 588 & 184 & 405 & 666 & \textbf{177} & 49.2 & 84.3 & 49.6 & 80.2 & \textbf{22.0} & 481 \\
5040 & 17920 & 560 & 580 & 178 & 400 & 656 & \textbf{176} & 52.7 & 88.0 & 53.2 & 86.8 & \textbf{23.1} & 561 \\
5760 & 20480 & 640 & 590 & 189 & 406 & 667 & \textbf{177} & 73.1 & 125.8 & 67.8 & 109.9 & \textbf{29.9} & 644 \\
6480 & 23040 & 720 & 595 & 216 & 410 & 673 & \textbf{178} & 80.2 & 136.9 & 79.1 & 130.8 & \textbf{34.2} & 725 \\
7200 & 25600 & 800 & 590 & 192 & 407 & 668 & \textbf{177} & 97.2 & 181.3 & 98.5 & 160.0 & \textbf{42.6} & 802 \\ \hline
\end{tabular}
\end{sidewaystable}

We further conduct validation on binary classification tasks using LIBSVM datasets \cite{Chang2011LIBSVMAL} converted to the MATLAB format. Due to the classification task, we reduced the penalty coefficient $\mu$ to $5 \times 10^{-4}$. Here, experiments are initialized at \(x^{0}=0\) with convergence criteria defined as:

\[
\frac{\|x^{n}-x^{n-1}\|}{\max\{1,\|x^{n}\|\}}<\epsilon\text{, where }\epsilon= 10^{-i}\text{, for }i=4,5,\ldots,9.
\]

\begin{sidewaystable}[htbp!]
\centering
\caption{Solving (\ref{scad}) on LIBSVM datasets\\ (\textcircled{1}DCA \textcircled{2}BDCA \textcircled{3}BapDCA \textcircled{4}3BapDCA \textcircled{5}3BapDCA$_\text{e}$)}
\label{table3}
\begin{tabular}{c|rrrrr|rrrrr}
\hline
\multicolumn{1}{c|}{Precision} & \multicolumn{5}{c|}{iter} & \multicolumn{5}{c}{CPU time (s)} \\ \hline
$\epsilon$ &\textcircled{1} & \textcircled{2} & \textcircled{3} & \textcircled{4} & \textcircled{5} &\textcircled{1} & \textcircled{2} & \textcircled{3} & \textcircled{4} & \textcircled{5} \\ \hline
$10^{-4}$ & 1216 & 1056 & 1583 & 1214 & \textbf{202} & 10.8 & 51.1 & 18.6 & 14.3 & \textbf{2.4} \\
$10^{-5}$ & 20644 & 7705 & 20799 & 20629 & \textbf{1003} & 210.2 & 351.2 & 281.2 & 262.0 & \textbf{11.6} \\
$10^{-6}$ & 168124 & 47693 & 175020 & 168305 & \textbf{10002} & 1190.7 & 1177.0 & 1531.5 & 1472.7 & \textbf{88.0} \\
$10^{-7}$ & Max & Max & Max & Max & \textbf{44403} & —— & —— & —— & —— & \textbf{471.0} \\
$10^{-8}$ & Max & Max & Max & Max & \textbf{250003} & —— & —— & —— & —— & \textbf{3415.1} \\
$10^{-9}$ & Max & Max & Max & Max & \textbf{410003} & —— & —— & —— & —— & \textbf{3897.7} \\
\hline
\end{tabular}
\end{sidewaystable}

Table \ref{table3} provides a numerical comparison of various algorithms for solving the SCAD model (3.1) using LIBSVM datasets. The results indicate that the proposed 3BapDCA$_{e}$ algorithm exhibits significantly superior computational efficiency when compared to alternative methods.

\subsection{Graphic Ginzburg-Landau model}
We now focus on a segmentation problem with graphic Ginzburg-Landau modeling, which deviates from the traditional phase-field model by integrating a prior term, thereby facilitating a semi-supervised assignment. The problem is articulated as follows \cite{BF}:
\begin{equation}\label{eq:nonlocalGL}
    \min_{x\in \mathbb{R}^N}E(x) = \sum_{i,j}\frac{\epsilon}{2}w_{ij}(x(i)-x(j))^2 +\frac{1}{\epsilon}\mathbb{W}(x) + \frac{\eta}{2}\sum_i\Lambda(i)(x(i)-y(i))^2
\end{equation}
where $x$ represents the image or data, indexed by \( i \) and \( j \). The energy functional combines a double-well potential \( \mathbb{W}(x) = \frac{1}{4} \sum_{i=1}^N (x(i)^2 - 1)^2 \), a prior term $\Lambda(i)$ weighted by \( \eta \), and a nonlocal interaction term with parameter \( \epsilon \). The weights \( w_{ij} = K(i,j) \cdot N(i,j) \) is determined by feature similarity \( K(i,j) = \exp(-\|P_i - P_j\|_2^2/\sigma^2) \) and proximity \( N(i,j) \), where $P_{i}$ serving as the feature of the data point $i$ and $\sigma^2$ controlling the kernel width. The matrix \( \Lambda \) and vector \( y \) encode prior knowledge. See \cite{Shensun2023} for the graph Laplacian construction.

For image segmentation, we set the model parameters as $\epsilon = \eta = 10$. We employ Algorithm \ref{alg:cap} to solve the nonconvex minimization problem with energy functional:

\begin{equation*}
E(x) = H(x) + F(x)
\end{equation*}
where
\begin{align*}
H(x) &= \sum_{ij}\frac{\epsilon}{2}w_{ij}(x(i)-x(j))^2 + \frac{\eta}{2}\sum_i\Lambda(i)(x(i)-y(i))^2, \quad 
F(x) = \frac{1}{\epsilon}\mathbb{W}(x).
\end{align*}

At each iteration $n$, we implement an implicit-explicit scheme by reformulating the energy as $H^n(x) - F^n(x)$. The time step $\delta t$ must satisfy $\delta t < \frac{8}{77L}$, with $L$ being the Lipschitz constant of $F(x)$. For other algorithms in Table \ref{tab:nonlocalGLmodel}, we use the following splitting:

\begin{equation*}
\begin{aligned}
    E(u) = &\sum_{i,j}\frac{\epsilon}{2}w_{ij}(u(i)-u(j))^2 + \frac{\eta}{2}\sum_i\Lambda(i)(u(i)-y(i))^2 + \frac{L}{2}\sum_i u(i)^2 \\
    &- \left(\frac{L}{2}\sum_i u(i)^2 +\frac{1}{\epsilon}\mathbb{W}(u)\right).
\end{aligned}
\end{equation*}

Table \ref{tab:nonlocalGLmodel} provides a comparison of 5 algorithms: DCA \cite{LeThi2018},  BDCA \cite[Section 3]{Artacho2018}, BapDCA \cite[Section 3]{Artacho2018}, 3BapDCA and our Algorithm, 3BapDCA$_{\text{e}}$. Specifically, algorithms BapDCA, 3BapDCA (no extrapolation), and 3BapDCA$_{\text{e}}$ employ a preconditioned scheme that utilizes only 10 iterations of a perturbed Jacobi preconditioner, which is a parallel-friendly approach well-suited for modern computing architectures. The first criterion in Table \ref{tab:nonlocalGLmodel}, DICE Bound, utilizes the DICE similarity coefficient to assess segmentation quality and is computed as:

\begin{equation}
\text{DICE} = \frac{2|X \cap Y|}{|X| + |Y|},
\end{equation}
where $|X|$ and $|Y|$ represent the pixel counts of the segmentation result and ground truth, respectively. Our DICE Bound here refers to the DICE coefficient reaching a value of 0.98, representing a high agreement between the segmentation and reference data.

The results presented in Table \ref{tab:nonlocalGLmodel} clearly demonstrate the promising performance of 3BapDCA$_{\text{e}}$, as it achieves the fastest convergence and lowest computational time compared to competing algorithms when meeting the stopping conditions. In some low-accuracy scenarios, the BDCA algorithm exhibits the fewest iterations among all algorithms, due to the high accuracy of its subproblem solution and the line search. However, the running time of 3BapDCA${\text{e}}$ remains the shortest. This demonstrates the high efficiency of 3BapDCA$_{\text{e}}$ across different precision requirements.

\begin{table}[htbp!]
\centering
\caption{Solving \ref{eq:nonlocalGL} on 11 different termination criteria. (Criteria \uppercase\expandafter{\romannumeral1}: $\|\nabla E(u)\|$, Criteria \uppercase\expandafter{\romannumeral2}: $\|x^n - x^{n-1}\|$)\\ \textcircled{1}DCA \textcircled{2}BDCA \textcircled{3}BapDCA \textcircled{4}3BapDCA \textcircled{5}3BapDCA$_{\text{e}}$}\label{tab:nonlocalGLmodel}
\begin{tabular}{cccrrrrr}
\hline
\multicolumn{3}{c}{Criteria} & \textcircled{1} & \textcircled{2} & \textcircled{3} & \textcircled{4} & \textcircled{5} \\ \hline
\multicolumn{2}{c}{\multirow{2}{*}{DICE Bound}} & Iter & 67 & 160 & 132 & 182 & \textbf{52} \\ 
& & Time(s) & 17.86 & 12.08 & 16.87 & 23.30 & \textbf{9.63} \\ \hline
\multirow{10}{*}{\uppercase\expandafter{\romannumeral1}} & \multirow{2}{*}{$10^{-1}$} & Iter & 249 & 90 & 208 & 286 & \textbf{81} \\
& & Time(s) & 26.56 & 14.51 & 25.12 & 35.30 & \textbf{12.95} \\ 
& \multirow{2}{*}{$10^{-2}$} & Iter & 597 & \textbf{182} & 518 & 717 & 190 \\
& & Time(s) & 59.57 & 25.86 & 58.30 & 84.09 & \textbf{25.35} \\ 
& \multirow{2}{*}{$10^{-3}$} & Iter & 1084 & 336 & 883 & 1223 & \textbf{318} \\
& & Time(s) & 107.17 & 49.04 & 97.35 & 141.34 & \textbf{39.92} \\ 
& \multirow{2}{*}{$10^{-4}$} & Iter & 3948 & 1218 & 3199 & 4429 & \textbf{908} \\
& & Time(s) & 380.75 & 202.80 & 345.21 & 504.24 & \textbf{107.08} \\ 
& \multirow{2}{*}{$10^{-5}$} & Iter & 5777 & 2428 & 3690 & 5105 & \textbf{1059} \\
& & Time(s) & 566.42 & 399.79 & 397.86 & 580.89 & \textbf{124.27} \\ \hline
\multirow{10}{*}{\uppercase\expandafter{\romannumeral2}} & \multirow{2}{*}{$10^{-1}$} & Iter & 123 & 88 & 114 & 137 & \textbf{83} \\
& & Time(s) & 14.57 & 14.25 & 15.07 & 18.43 & \textbf{13.18} \\ 
& \multirow{2}{*}{$10^{-2}$} & Iter & 380 & \textbf{138} & 470 & 422 & 172 \\
& & Time(s) & 39.93 & 25.62 & 53.17 & 50.69 & \textbf{23.31} \\ 
& \multirow{2}{*}{$10^{-3}$} & Iter & 1044 & 328 & 861 & 1169 & \textbf{310} \\
& & Time(s) & 101.98 & 47.76 & 94.99 & 135.24 & \textbf{39.01} \\ 
& \multirow{2}{*}{$10^{-4}$} & Iter & 3749 & 1176 & 3087 & 4166 & \textbf{402} \\
& & Time(s) & 361.74 & 195.56 & 333.23 & 474.48 & \textbf{49.48} \\ 
& \multirow{2}{*}{$10^{-5}$} & Iter & 4363 & 1552 & 3579 & 4861 & \textbf{1002} \\
& & Time(s) & 420.39 & 255.64 & 385.93 & 553.20 & \textbf{117.78} \\ \hline
\end{tabular}
\end{table}

\begin{figure}[htbp]   
  \centering            
  \subfloat[redflower]   
  {\label{fig:redflower}\includegraphics[width=0.18\textwidth]{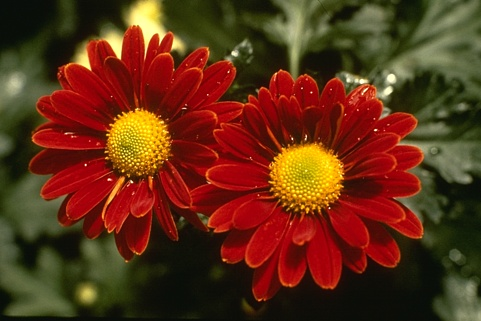}}\quad
  \subfloat[Label]
  {\label{fig:redflowerlabel}\includegraphics[width=0.18\textwidth]{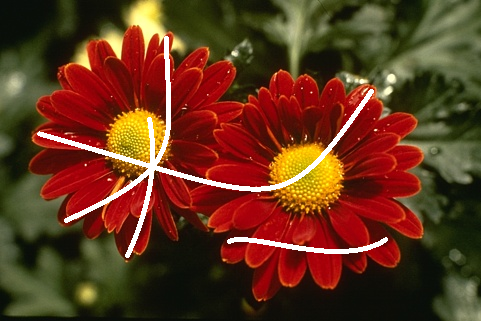}}\quad
  \subfloat[result 1]
  {\label{fig:result1}\includegraphics[width=0.18\textwidth]{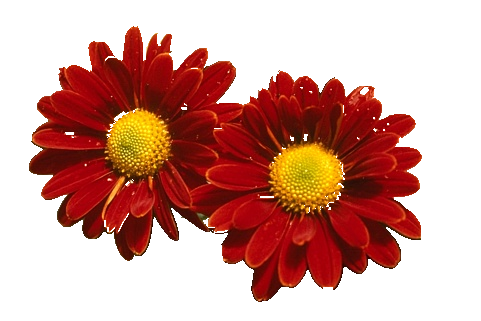}}\quad
  \subfloat[result 2]
  {\label{fig:result2}\includegraphics[width=0.18\textwidth]{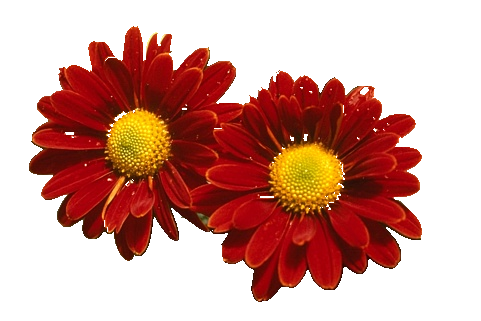}}\quad
  \subfloat[result 3]
  {\label{fig:result3}\includegraphics[width=0.18\textwidth]{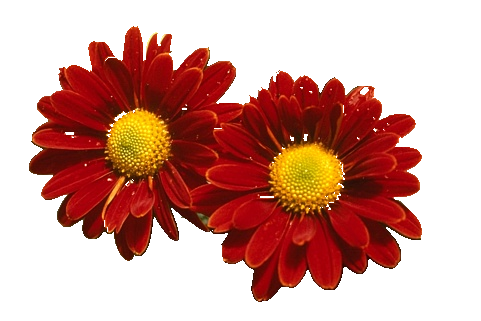}}\quad
  \caption{The performance of segmentation assignment.}
  \label{fig:redflower_seg}         
\end{figure}

In Figure \ref{fig:redflower_seg}, we illustrate the segmentation result of an image. Figure \eqref{fig:redflower} is the original image. Figure \eqref{fig:redflowerlabel} labels the only prior of the red flower to be segmented with the white pixels label and $\Lambda(i)=1$ (otherwise $\Lambda(i)=0$ and here we do not need the prior of the background). The last three figures represent the segmentation results under three different criteria (DICE Bound, $\|\nabla E(u)\|\leq 10^{-5}$ and $\|x^n-x^{n-1}\|\leq 10^{-5}$).

We finally end this section with a remark for the KL properties of the modified SCAD regularization \eqref{scad} and the graphic Ginzburg-Landau functional \eqref{eq:nonlocalGL}.

\begin{remark}\label{KL:two:models}
The discrete graph Ginzburg-Landau functional, being a polynomial in $x$, is semi-algebraic and consequently satisfies the Kurdyka-Łojasiewicz (KL) property \cite[Section 2.2]{ABS}. For the modified SCAD regularization \eqref{eq:huber-scad}, each component $p_{M,i}(x_i)$ constitutes a one-dimensional piecewise quadratic function. Following analogous reasoning to \cite[Section 5.2]{Li2018}, we establish that the energy functional $E(x)$ is a KL function. Moreover, application of \cite[Theorem 3.6]{Li2018} demonstrates that $A(x,y,z)$ is a KL function. These properties collectively guarantee the convergence for both models under consideration.
\end{remark}

\section{Conclusion}\label{sec:conclusion}
We propose a preconditioned third-order implicit-explicit algorithm with extrapolation (3BapDCA$_{\text{e}}$) for solving \eqref{eq:basefunctional}. Our algorithmic framework allows different choices of extrapolation parameters $\{\beta_n\}_n$. We establish the global and local convergence of the consequence generated by 3BapDCA$_{\text{e}}$ by assuming the Kurdyka-\L ojasiewicz property of the objective and the L-smoothness of $F$. Our numerical experiments verify that our algorithm outperforms other classic algorithms for least squares problems with the SCAD regularizer and the graphic Ginzburg-Landau model.

\noindent
{\small
	\textbf{Acknowledgements}
 Kelin Wu and Hongpeng Sun acknowledge the support of the National Key R\&D Program of China (2022ZD0116800),  the National Natural Science Foundation of China under grant No. \,12271521,  and the Beijing Natural Science Foundation No. Z210001. 
}

\bibliographystyle{plain}
\bibliography{bdfab_ls}
\end{document}